\documentclass{conm-p-l}
\usepackage{amssymb}
\usepackage[mathscr]{eucal}
\usepackage{amscd}
\usepackage{pb-diagram}

\def\Bbb{\mathbb}
\def\frak{\mathfrak}

\newenvironment{pf*}[1]{\proof[#1]}{\endproof}
\newtheorem{Theorem}[equation]{Theorem}

\newtheorem{Lemma}[equation]{Lemma}
\newtheorem{Proposition}[equation]{Proposition}

\theoremstyle{definition}
\newtheorem{Definition}[equation]{Definition}
\newtheorem{Example}[equation]{Example}

\theoremstyle{remark}
\newtheorem{Remark}[equation]{Remark}

\numberwithin{equation}{section}

\newcommand{\thmref}[1]{Theorem~\ref{#1}}
\newcommand{\secref}[1]{\S\ref{#1}}
\newcommand{\lemref}[1]{Lemma~\ref{#1}}
\newcommand{\propref}[1]{Proposition~\ref{#1}}

\newcommand{\defref}[1]{Definition~\ref{#1}}

\newcommand{\defeq}{\overset{\operatorname{\scriptstyle def.}}{=}}
\newcommand{\C}{{\Bbb C}}
\newcommand{\Z}{{\Bbb Z}}
\newcommand{\Q}{{\Bbb Q}}

\newcommand{\g}{{\frak g}}

\newcommand{\ve}{q}
\newcommand{\vaep}{\varepsilon}
\newcommand{\M}{{\frak M}} 

\newcommand{\La}{{\frak L}} 
\newcommand{\bv}{{\mathbf v}} 
\newcommand{\bw}{{\mathbf w}} 

\newcommand{\Uq}{{\mathbf U}_q(\mathfrak g)} 

\newcommand{\Ua}{{\mathbf U}_q(\widehat{\mathfrak g})} 
\newcommand{\Ul}{{\mathbf U}_q({\mathbf L}{\mathfrak g})} 

\newcommand{\Lg}{\mathbf L\g}

\newcommand{\te}{{\widetilde e}}
\newcommand{\tf}{{\widetilde f}}
\newcommand{\wt}{\operatorname{wt}}
\newcommand{\B}{\mathcal B}
\newcommand{\Y}[2]{Y_{#1,#2}}

\newcommand{\bfR}{\mathbf R}

\newcommand{\qch}[1]{\widetilde{\chi_{#1}}}

\newcommand{\yl}{14pt}
\newcommand{\yh}{7pt}
\newcommand{\ffbox}[1]{
\setbox9=\hbox{$\scriptstyle\overline{1}$}
\framebox[\yl][c]{\rule{0mm}{\ht9}${\scriptstyle #1}$}
}
\newcommand{\numbullet}[1]{
\setbox8=\hbox{${#1}$}
\setbox9=\hbox{$\overset{#1}\bullet$}
\raisebox{\ht8}{\raisebox{-\ht9}{\box9}}
}
\newcommand{\fhbox}[1]{
\setbox9=\hbox{$\scriptstyle\overline{1}$}
\framebox[\yh][c]{\rule{0mm}{\ht9}${\scriptstyle #1}$}
}
\newcommand{\twoY}[2]{
\renewcommand{\arraycolsep}{0pt}
\begin{array}{|c|}
   \hline \hbox to \yl{\hfill$\scriptstyle #1$\hfill} \\
   \hline \hbox to \yl{\hfill$\scriptstyle #2$\hfill} \\
   \hline
\end{array}
}

\begin{document}
\title[quantum affine algebras of type $A_n$, $D_n$]
{$t$--analogs of $q$--characters of quantum affine algebras of type
$A_n$, $D_n$
}
\author{Hiraku Nakajima}
\address{Department of Mathematics, Kyoto University, Kyoto 606-8502,
Japan
}
\email{nakajima@kusm.kyoto-u.ac.jp}
\urladdr{http://www.kusm.kyoto-u.ac.jp/\textasciitilde nakajima}
\thanks{Supported by the Grant-in-aid
for Scientific Research (No.13640019), JSPS}
%
\dedicatory{Dedicated to Professor Ryoshi Hotta on his sixtieth birthday}
\copyrightinfo{2002}
    {American Mathematical Society}
\begin{abstract}
We give a tableaux sum expression of $t$--analog of $q$--characters of
finite dimensional representations (standard modules) of quantum
affine algebras $\Ul$ when $\g$ is of type $A_n$, $D_n$.
\end{abstract}
\maketitle

\section{Introduction}

Let $\g$ be a simple Lie algebra of type $ADE$ over $\C$, $\Lg = \g
\otimes \C[z,z^{-1}]$ be its loop algebra, and $\Ul$ be its quantum
universal enveloping algebra, or the quantum loop algebra for short.
It is a subquotient of the quantum affine algebra $\Ua$, i.e., without
central extension and degree operator.
It is customary to define $\Ul$ as an algebra over $\Q(q)$, but here
we consider $q$ as a nonzero complex number which is not a root of
unity, for simplicity.

By Drinfeld \cite{Dr}, Chari-Pressley \cite{CP-rep}, simple
$\Ul$-modules are parametrized by $I$-tuples of polynomials $P =
(P_i(u))_{i\in I}$ with normalization $P_i(0) = 1$. They are called
{\it Drinfeld polynomials}. Let us denote by $L(P)$ the simple module
with Drinfeld polynomial $P$. It gives a basis $\{ L(P) \}_P$ of the
Grothendiek group $\operatorname{Rep}\Ul$ of the category of finite
dimensional representations of $\Ul$.

In \cite{Na-qaff} the author introduced another set of $\Ul$-modules
$M(P)$, called {\it standard modules\/}, parametrized also by Drinfeld
polynomials. It gives us another base of $\operatorname{Rep}\Ul$.
Then the author \cite{Na-qaff} showed that the multiplicity
$[M(P):L(Q)]$ is equal to a specialization of a polynomial
$Z_{PQ}(t)\in\Z[t,t^{-1}]$ at $t=1$. And the polynomial $Z_{PQ}(t)$ is
defined as Poincar\'e polynomials of intersection cohomology of graded
quiver varieties, which are fixed point sets of $\C^*$-actions on
quiver varieties, introduced earlier by the author
\cite{Na:1994,Na:1998}.

The polynomials $Z_{PQ}(t)$ can be considered as an analog of
Kazhdan-Lusztig polynomials which are Poincar\'e polynomials of
intersection cohomology of Schubert varieties. As Kazhdan-Lusztig
polynomials are defined via an involution on the Hecke algebra, our
$Z_{PQ}(t)$ are determined by means of a bar involution
$\setbox5=\hbox{A}\overline{\rule{0mm}{\ht5}\hspace*{\wd5}}\,$ on the
$t$--analog of the Grothendieck ring
$\bfR_t \defeq \operatorname{Rep}\Ul\otimes_\Z\Z[t, t^{-1}]$.
In order to compute the bar involution (and hence $Z_{PQ}(t)$), the
author introduced $t$--analogs of $q$--characters $\qch{q,t}$ and
gave a combinatorial algorithm to compute $\qch{q,t}(M(P))$
\cite{Na-qchar,Na-qchar-main}. 
The original $q$--characters $\chi_{q}$ had been introduced and studied
by Knight, Frenkel-Reshetikhin, Frenkel-Mukhin \cite{Kn,FR,FM}.
In summary, the multiplicity $[M(P):L(Q)]$ can be given by a purely
combinatorial algorithm.

In this paper, we shall give an {\it explicit\/} expresseion of
$\qch{q,t}(M(P))$ when $\g$ is of type $A_n$, $D_n$ in terms of Young
tableaux or their variants. Such expressions had been known for the
original $q$--characters $\chi_{q}$ by Kuniba-Suzuki \cite{KuS}
and Nazarov-Tarasov \cite{NT}.
(They did not use the terminology of $q$--characters. But their
calculation can be translated to $q$--characters. See \cite{FR1},
\cite[\S11]{FR} and reference therein.)
In fact, the author finds this expression via a certain relation
between $q$--characters and Kashiwara's crystal base (see
\secref{sec:crystal}), where expression in terms of tableaux was given
by Kashiwara-Nakashima \cite{Kas-Na}. (See also \cite{Kang}.) It seems that
the relation between crystals and $\chi_q$ has not been known. 
The author is also motivated by his earlier work \cite{Na:Hom} on an
expression of Poincar\'e polynomials of original quiver varieties of
type $A_n$ in terms of Young tableaux. This work was motivated by
works of Shimomura, Hotta-Shimomura \cite{S,HS} in turn.

In this paper we use the following notation: $(P)$ be $1$ if a
statement $P$ is true and $0$ otherwise.

\section{$t$--analogs of $q$--characters}

We shall not discuss the definition of quantum loop algebras, nor
their finite dimensional representations in this paper. (See
\cite{Na-qchar} for a survey.) We just review properties of
$\qch{q,t}$, as axiomized in \cite{Na-qchar-main}.

Let $\g$ be a simple Lie algebra of type $ADE$, let $I$ be the index
set of simple roots.
Let $L(P)$ (resp.\ $M(P)$) be the simple (resp.\ standard)
$\Ul$-module with Drinfeld polynomial $P$.
A simple module $L(P)$ is called an {\it l--fundamental
representation\/} when $P_i(u) = (1 - au)^{\delta_{iN}}$ for some
$a\in\C^*$ and $N\in I$. Since it depends only on $N$ and $a$, we
denote it by $L(\Lambda_N)_a$. This will play an important role later.

Let
\(
\mathscr Y_t \defeq 
   \Z[t,t^{-1},Y_{i,a}, Y_{i,a}^{-1}]_{i\in I, a\in\C^*}
\)
be a Laurent polynomial ring of uncontably many variables $Y_{i,a}$'s
with coefficients in $\Z[t,t^{-1}]$. A {\it monomial\/} in $\mathscr
Y_t$ means a monomial only in $Y_{i,a}^\pm$, containing no $t$'s.
Let
\begin{equation*}
   A_{i,a} \defeq Y_{i,a\ve} Y_{i,a\ve^{-1}}
     \prod_{j:j\neq i} Y_{j,a}^{c_{ij}},
\end{equation*}
where $c_{ij}$ is the $(i,j)$-entry of the Cartan matrix.
Let $\mathcal M$ be the set of monomials in $\mathscr Y_t$.

\begin{Definition}
(1) For a monomial $m\in\mathcal M$, we define $u_{i,a}(m)\in\Z$ be the
degree in $Y_{i,a}$, i.e.,
\begin{equation*}
   m = \prod_{i,a} Y_{i,a}^{u_{i,a}(m)}.
\end{equation*}

(2) A monomial $m\in\mathcal M$ is said $i$--dominant if
$u_{i,a}(m)\ge 0$ for all $a$. It is said {\it l--dominant} if
it is $i$--dominant for all $i$.

(3) Let $m, m'$ be monomials in $\mathcal M$. We say $m \le m'$ if
$m/m'$ is a monomial in $A_{i,a}^{-1}$ ($i\in I$, $a\in\C^*$).
Here a monomial in $A_{i,a}^{-1}$ means a product of nonnegative
powers of $A_{i,a}^{-1}$. It does not contain any factors
$A_{i,a}$. In such a case we define $v_{i,a}(m, m')\in\Z_{\ge 0}$ by
\begin{equation*}
   m = m' \prod_{i,a} A_{i,a}^{-v_{i,a}(m,m')}.
\end{equation*}
This is well-defined since the $\ve$-analog of the Cartan matrix is
invertible. We say $m < m'$ if $m\le m'$ and $m\neq m'$.

(4) For an $i$--dominant monomial $m\in\mathcal M$ we define
\begin{equation*}
   E_i(m) \defeq
    m\, \prod_a
     \sum_{r_a=0}^{u_{i,a}(m)}
     t^{r_a(u_{i,a}(m)-r_a)}
     \begin{bmatrix}
       u_{i,a}(m) \\ r_a
     \end{bmatrix}_t A_{i,a\ve}^{-r_a},
\end{equation*}
where
\(
\left[\begin{smallmatrix}
  n \\ r
\end{smallmatrix}\right]_t
\)
is the $t$-binomial coefficient.
\end{Definition}

Suppose that {\it l\/}--dominant monomials $m_{P^1}$, $m_{P^2}$ and
monomials $m^1\le m_{P^1}$, $m^2\le m_{P^2}$ are given. We define an
integer $d(m^1, m_{P^1}; m^2, m_{P^2})$ by
\begin{multline}\label{eq:d}
   d(m^1, m_{P^1}; m^2, m_{P^2})
\\
   \defeq
   \sum_{i,a} \left( v_{i,a\ve}(m^1, m_{P^1}) u_{i,a}(m^2)
   + u_{i,a\ve}(m_{P^1}) v_{i,a}(m^2, m_{P^2})\right).
\end{multline}

For an $I$-tuple of rational functions $Q/R = (Q_i(u)/R_i(u))_{i\in
I}$ with $Q_i(0) = R_i(0) = 1$, we set
\begin{equation*}
   m_{Q/R} \defeq
   \prod_{i\in I} \prod_{\alpha} \prod_{\beta}
     Y_{i,\alpha} Y_{i,\beta}^{-1},  
\end{equation*}
where $\alpha$ (resp.\ $\beta$) runs roots of $Q_i(1/u) = 0$
(resp.\ $R_i(1/u) = 0$), i.e.,
$Q_i(u) = \prod_\alpha ( 1 - \alpha u)$ (resp.\ $R_i(u) = \prod_\beta
(1 - \beta u)$). As a special case, an $I$-tuple of polynomials $P =
(P_i(u))_{i\in I}$ defines $m_P = m_{P/1}$.
In this way, the set $\mathcal M$ of monomials are identified with the 
set of $I$-tuple of rational functions, and the set of {\it
l\/}--dominant monomials are identified with the set of $I$-tuple of
polynomials.

The $t$--analog of the Grothendieck ring $\bfR_t$ 
is a free $\Z[t,t^{-1}]$-module with base $\{ M(P) \}$ where $P =
(P_i(u))_{i\in I}$ is the Drinfeld polynomial.
(We do not recall the definition of standard modules $M(P)$ here, but
the reader safely consider them as formal variables.)

The $t$--analog of the $\ve$--character homomorphism is a
$\Z[t,t^{-1}]$-linear homomorphism
\(
   \qch{q,t}\colon \bfR_t \to \mathscr Y_t.
\)
It is defined as the generating function of Poincar\'e polynomials of
graded quiver varieties, and will not be reviewed in this paper.
But the following is known.

\begin{Theorem}\label{thm:ind}
\textup{(1)} 
The $\qch{q,t}$ of a standard module $M(P)$ has a form
\begin{equation*}
   \qch{q,t}(M(P)) = m_P + \sum a_m(t) m,
\end{equation*}
where the summation runs over monomials $m < m_P$.

\textup{(2)}
For each $i\in I$, $\qch{q,t}(M(P))$ can be expressed as a linear
combination \textup(over $\Z[t,t^{-1}]$\textup) of $E_i(m)$ with
$i$--dominant monomials $m$.

\textup{(3)}
Suppose that two $I$-tuples of polynomials $P^1 = (P^1_i)$, $P^2 =
(P^2_i)$ satisfy the following condition:
\begin{equation}
\label{eq:Z}
\begin{minipage}[m]{0.75\textwidth}
\noindent   
   $a/b\notin\{ \ve^n \mid n\in\Z, n \ge 2\}$ for any
   pair $a$, $b$ with $P^1_i(1/a) = 0$, $P^2_j(1/b) =
   0$ \textup($i,j\in I$\textup).
\end{minipage}
\end{equation}
Then we have
\begin{equation*}
  \qch{q,t}(M(P^1P^2)) =
  \sum_{m^1, m^2} t^{2d(m^1, m_{P^1}; m^2, m_{P^2})}
   a_{m^1}(t) a_{m^2}(t) m^1 m^2
,
\end{equation*}
where
\(
   \qch{q,t}(M(P^a)) = \sum_{m^a} a_{m^a}(t) m^a
\)
with $a=1,2$.

Moreover, properties \textup{(1),(2),(3)} uniquely determine $\qch{q,t}$.
\end{Theorem}

Apart from the existence problem, one can consider the above
properties (1), (2), (3) as the definition of $\qch{q,t}$ (an
axiomatic definition). We only use the above properties, and the
reader can safely forget the original definition.

\begin{Remark}
  It is more suitable to consider a slightly modified version
  $\chi_{q,t}$ in stead of $\qch{q,t}$ for computing the bar
  operation. Therefore $\chi_{q,t}$ was mainly used in \cite{Na-qchar}.
  Anyhow, they are simply related as
\[
  \text{if\ } \qch{q,t}(M(P)) = \sum_m a_m(t) m , \qquad
  \chi_{q,t}(M(P)) = \sum_m t^{-d(m,m_P;m, m_P)} a_m(t) m
  .
\]
See \cite[5.1.3]{Na-qchar}.
\end{Remark}

Let us explain briefly why the properties (1), (2), (3) determine
$\qch{q,t}$. First consider the case $M(P)$ is an {\it
  l\/}--fundamental representation. (We have $M(P) = L(P)$ in this
case.) Then one can determine $\qch{q,t}(M(P))$ starting from $m_P$
and using the property (2) inductively. (The idea can be seen in the
examples below.)
For general $P$, write it as $P = P^1 P^2 P^3 \cdots$ so that each
$M(P^\alpha)$ is an {\it l\/}--fundamental representation, and
the condition~\eqref{eq:Z} is met with respect to the ordering. Then
we apply (3) successively to get
\(
   \qch{q,t}(M(P))
\)
from
\(
   \qch{q,t}(M(P^\alpha))
\)
with $\alpha=1,2,\cdots$.

We attach to each standard module $M(P)$, an oriented colored graph
$\Gamma_P$ as follows. (It is a slight modification of the graph in
\cite[5.3]{FR}.) The vertices are monomials in $\qch{q,t}(M(P))$. We
draw a colored edge $\xrightarrow{i,a}$ from $m_1$ to $m_2$ if $m_2 =
m_1 A_{i,a}^{-1}$. We also write the coefficients of the monomials in
$\qch{q,t}(M(P))$. In fact, edges are determined from monomials on
vertices.

Here are examples.
\begin{Example}\label{ex:graph}
Let $\g$ be of type $A_2$. We put a numbering $I = \{ 1, 2\}$.

(1) The graph of $\qch{q,t}(M(P))$ with $m_P = \Y 11 \Y 2\ve$ is the
following:
\begin{equation*}
\begin{CD}
  \Y 11 \Y 2\ve @>{2,\ve^2}>> \Y 11 \Y 1{\ve^2} \Y2{\ve^3}^{-1}
  @>{1,\ve^3}>> \Y 11 \Y 1{\ve^4}^{-1} @>{1,\ve}>>
  \Y 1{\ve^2}^{-1}\Y 1{\ve^4}^{-1}\Y 2\ve \\
  @V{1,\ve}VV @V{1,\ve}VV @. @V{2,\ve^2}VV \\
  \Y 1{\ve^2}^{-1}\Y 2{\ve}^2 @>{2,\ve^2}>>
  (1 + t^2) \Y 2{\ve} \Y 2{\ve^3}^{-1} @>{2,\ve^2}>>
  \Y 1{\ve^2} \Y 2{\ve^3}^{-2}
  @>{1,\ve^3}>> \Y 1{\ve^4}^{-1}\Y 2{\ve^3}^{-1}
\end{CD}
\end{equation*}

(2) The graph of $\qch{q,t}(M(P))$ with $m_P = \Y 11^2$ is the
following:
\begin{equation*}
\begin{CD}
  \Y 11^2 @>{1,\ve}>> (1+t^2)\Y 11 \Y 1{\ve^2}^{-1} \Y2\ve
  @>{1,\ve}>> \Y 1{\ve^2}^{-2} \Y 2\ve^2 @. \\
  @. @V{2,\ve^2}VV @V{2,\ve^2}VV \\
  @. (1+t^2)\Y 11\Y 2{\ve^3}^{-1} @>{1,\ve}>>
  (1+t^2)\Y 1{\ve^2}^{-1}\Y 2\ve \Y 2{\ve^3}^{-1}
  @>{2,\ve^2}>> \Y 2{\ve^3}^{-2}
\end{CD}
\end{equation*}

(3) The graph of $\qch{q,t}(M(P))$ with $m_P = \Y 11 \Y 1{\ve^2}$ is
the following:
\begin{equation*}
\begin{CD}
  \Y 11 \Y 1{\ve^2} @>{1,\ve}>> \Y 2\ve @>{2,\ve^2}>>
    \Y 1{\ve^2} \Y 2{\ve^3}^{-1}
\\
  @V{1,\ve^3}VV @. @V{1,\ve^3}VV
\\
  \Y 11 \Y 1{\ve^4}^{-1} \Y 2{\ve^3} @>{1,\ve}>>
  \Y 1{\ve^2}^{-1} \Y 1{\ve^4}^{-1} \Y 2{\ve} \Y 2{\ve^3}
  @>{2,\ve^2}>> \Y 1{\ve^4}^{-1}
\\
  @V{2,\ve^4}VV @V{2,\ve^4}VV @.
\\
  \Y 11\Y 2{\ve^5}^{-1} @>{1,\ve}>>
  \Y 1{\ve^2}^{-1}\Y 2{\ve}\Y 2{\ve^5}^{-1}
  @>{2,\ve^2}>> \Y 2{\ve^3}^{-1} \Y 2{\ve^5}^{-1}
\end{CD}
\end{equation*}
\end{Example}

In the first two examples, $\qch{q,t}(M(P))$ does not contain {\it
  l\/}--dominant monomials other than $m_P$. Therefore the graphs are
determined from \thmref{thm:ind}(1), (2), as we mentioned above. (It
is instructive to check that (2) holds in these examples.) Other
examples can be found in \cite{Na-qchar}. ({\bf Caution}: (1) In
[loc.\ cit.], $\chi_{q,t}$ in stead of $\qch{q,t}$ was used. (2) There
are mistakes in Example~5.3.3 in [loc.\ cit.].)

\section{A monomial realization of crystal bases}\label{sec:crystal}

In this section, we give a realization of crystal bases of highest
weight modules, called a {\it monomial realization}. We can avoid the
usage of this material in later sections, but it will give us a
natural motivation for our tableaux sum expression of
$\ve$--characters. Moreover, this section can be read independently
from the other sections. See also Kashiwara's article \cite{Kas-Korea}
in this volume.

In this section, $\mathfrak g$ is an arbitrary symmetrizable Kac-Moody
Lie algebra. Let $I$ be the index set of simple roots. Let $\{
\alpha_i \}_{i\in I}$, $\{ h_i \}_{i\in I}$ be the sets of simple
roots and simple coroots. Let $P$ be a weight lattice, and $P^+$ be
the set of dominant weights. We assume that there exists $\Lambda_i\in
P$ such that $\langle h_i, \Lambda_j\rangle = \delta_{ij}$
(fundamental weights).

We shall not recall here the notion of crystals. See e.g.,
\cite{Kas-Dem}.
We denote by $\B(\lambda)$ the crystal of the highest weight module
with highest weight $\lambda\in P^+$.

Let $\mathcal M^\circ$ be the set of monomials in $\mathscr Y =
\Z[Y_{i,a}^\pm]_{i\in I, a\in\C^*}$ such that $a$ is a power of $\ve$:
\begin{equation*}
   \mathcal M^\circ \defeq \left\{\left.
     m = \prod_{i,n} Y_{i,\ve^n}^{u_{i,\ve^n}(m)}\, \right|
   \text{$u_{i,\ve^n}(m)\in\Z$ is zero except for finitely many $(i,n)$}
   \right\}.
\end{equation*}
We set
\begin{gather*}
   \vaep_{i,n}(m) \defeq - \sum_{k:k\ge n} u_{i,\ve^k}(m),
\qquad
   \varphi_{i,n}(m) \defeq \sum_{k:k \le n} u_{i,\ve^k}(m),
\\
  \vaep_i(m) \defeq \max_n \vaep_{i,n}(m), \qquad
  \varphi_i(m) \defeq \max_n \varphi_{i,n}(m),
\\
   p_i(m) \defeq \max \{ n \mid \vaep_{i,n}(m) = \vaep_i(m) \},
\qquad
   q_i(m) \defeq \min \{ n \mid \varphi_{i,n}(m) = \varphi_i(m) \},   
\\
  \wt(m) \defeq \sum_{i,n} u_{i,\ve^n}(m) \Lambda_i.
\end{gather*}
If $n$ is sufficiently large, we have $\vaep_{i,n}(m) = 0$. 
If $n$ is sufficiently small, $\vaep_{i,n}(m)$ is a fixed integer
independent of $n$. Therefore $\vaep_i(m)$ is a nonnegative integer. If
$\vaep_i(m) = 0$, we understand $p_i(m) = \infty$. Similarly, we set
$q_i(m) = -\infty$ if $\varphi_i(m) = 0$.
We define operators $\te_i$, $\tf_i$ by
\begin{equation*}
\begin{split}
   & \te_i(m) \defeq 
   \begin{cases}
     m A_{i,\ve^{p_i(m)-1}}^{-1} & \text{if $\vaep_i(m) > 0$}, \\
     0 & \text{if $\vaep_i(m) = 0$},
   \end{cases}
\\
  & \tf_i(m) \defeq
  \begin{cases}
     m A_{i,\ve^{q_i(m)+1}} & \text{if $\varphi_i(m) > 0$}, \\
     0 & \text{if $\varphi_i(m) = 0$}.
  \end{cases}
\end{split}
\end{equation*}

Unfortunately this does not give us a crystal in general. Therefore we
need an extra assumption or modification. Here we assume that
$\mathfrak g$ is {\it without odd cycles}, i.e., there exists a
function $I\in i\mapsto a_i\in \{0, 1\}$ such that $a_i + a_j = 1$
whenever $a_{ij} < 0$ here.
This condition is satisfied when $\mathfrak g$ is finite-dimensional
(so enough for our present purpose), or of affine type other than
$A_{2n}^{(1)}$. (See \cite[\S8]{Na-tensor} or \cite{Kas-Korea} for other
modifications of the rule to have a crystal for arbitrary $\g$.)

Let
\begin{equation*}
   \mathcal M' \defeq
   \left\{ \left. m = \prod_{i,n} Y_{i,\ve^n}^{u_{i,\ve^n}(m)}
   \in \mathcal M^\circ\, \right| 
 \text{$u_{i,\ve^n}(m) = 0$ if $n\equiv a_i\bmod 2$}
 \right\}.
\end{equation*}
It is clear that $\mathcal M'$ is invariant under $\te_i$, $\tf_i$.

\begin{Theorem}
\textup{(1)} The set $\mathcal M'$ together with maps $\wt$, $\vaep_i$,
$\varphi_i$, $\te_i$, $\tf_i$ satisfies the axioms of a crystal in
the sense of \cite{Kas-Dem}.

\textup{(2)} The crystal generated by an {\it l\/}--dominant monomial
$m\in\mathcal M'$ \textup(i.e., $u_{i,\ve^n}(m)\ge 0$ for all
$i,n$\textup) is isomorphic to the crystal $\B(\wt(m))$ of the highest
weight module.
\end{Theorem}

{\bf Warning}: The crystals, we constructed, are for $\g$, not for
$\widehat\g$ or ${\mathbf L}\g$.

Here are examples.
\begin{Example}[Compare with Example~\ref{ex:graph}]
Let $\g$ be of type $A_2$.

(1) The crystal graph of $\mathcal M'$ starting from $\Y 11 \Y 2\ve$
is the following:
\begin{equation*}
\begin{CD}
  \Y 11 \Y 2\ve @>2>> \Y 11 \Y 1{\ve^2} \Y2{\ve^3}^{-1} @>1>> 
  \Y 11 \Y 1{\ve^4}^{-1} 
  @>1>> \Y 1{\ve^2}^{-1}\Y 1{\ve^4}^{-1}\Y 2\ve \\
  @V1VV @. @. @V2VV \\
  \Y 1{\ve^2}^{-1}\Y 2{\ve}^2 @>>2> \Y 2{\ve} \Y 2{\ve^3}^{-1}
  @>2>>
  \Y 1{\ve^2} \Y 2{\ve^3}^{-2}
  @>1>> \Y 1{\ve^4}^{-1}\Y 2{\ve^3}^{-1}
\end{CD}
\end{equation*}

(2) The crystal graph of $\mathcal M'$ starting from $\Y 11^2$ is the
following: 
\begin{equation*}
\begin{CD}
  \Y 11^2 @>1>> \Y 11 \Y 1{\ve^2}^{-1} \Y2\ve @>1>>
  \Y 1{\ve^2}^{-2} \Y 2\ve^2 @. \\
  @. @V2VV @V2VV \\
  @. \Y 11\Y 2{\ve^3}^{-1} @>1>> \Y 1{\ve^2}^{-1}\Y 2\ve \Y 2{\ve^3}^{-1}
  @>2>> \Y 2{\ve^3}^{-2}
\end{CD}
\end{equation*}

(3) The crystal graph of $\mathcal M'$ starting from $\Y 11 \Y
1{\ve^2}$ is the following:
\begin{equation*}
\begin{CD}
  \Y 11 \Y 1{\ve^2} @. @. @.
\\
  @V{1}VV @. @.
\\
  \Y 11 \Y 1{\ve^4}^{-1} \Y 2{\ve^3} @>{1}>>
  \Y 1{\ve^2}^{-1} \Y 1{\ve^4}^{-1} \Y 2{\ve} \Y 2{\ve^3}
  @. @.
\\
  @V{2}VV @V{2}VV @.
\\
  \Y 11\Y 2{\ve^5}^{-1} @>{1}>>
  \Y 1{\ve^2}^{-1}\Y 2{\ve}\Y 2{\ve^5}^{-1}
  @>{2}>> \Y 2{\ve^3}^{-1} \Y 2{\ve^5}^{-1}
\end{CD}
\end{equation*}
Note that (2) and (3) are different realization of the same crystal
$\B(2\Lambda_1)$. Comparing with Example~\ref{ex:graph}, we find that
vertices and edges are subsets of those of graphs for
$\qch{q,t}$. This is the case for any {\it l\/}--dominant monomials
in $\mathcal M'$, as we shall explain below.
\end{Example}

There are two proofs of this theorem. The author's proof is based on
\cite[8.6]{Na-tensor}, in particular depends on the theory of quiver
varieties. There is more direct proof due to Kashiwara
\cite{Kas-Korea}\footnote{In fact, when the author obtained
\cite[8.6]{Na-tensor}, he thought it a new result on crystals. But
afterwards, Kashiwara informed him that it follows directly from the
main result of \cite{Kas-Dem} together with the formula
$T_\lambda\otimes\B_i = T_{s_i\lambda}\otimes\B_i$ (see \cite{Kas-Dem}
for the notation). The last setence of the abstract in
\cite{Na-tensor} should be corrected as `This result is equivalent to
Kashiwara's combinatorial description given by his embedding
theorem'.}.
We explain the author's proof here since it is related closely to
$\ve$--characters, although we shall not explain quiver varieties.
(See e.g., \cite[\S8]{Na-qchar} for a summary of the theory of quiver
varieties.)

The $t$--analogs of $\ve$--characters of standard modules are the
generating function of Poincar\'e polynomials of graded quiver
varieties, which are fixed points of the quiver variety $\M$ with
respect to a $\C^*$-action. It is expressed schematically as
\begin{equation*}
   \qch{q,t}(M(P)) = \sum_{m} P_t(\M(m)) m,
\end{equation*}
where $\M(m)$ is the connected component corresponding to a monomial
$m$, and $P_t(\M(m))$ is its Poincar\'e polynomial. The choice of the
$\C^*$-action corresponds to a choice of the Drinfeld polynomial
$P$. The {\it l\/}--dominant monomial $m_P$ corresponding to $P$ is a
certain distinguished component, which is a single point (and hence
$P_t(\M(m_P)) = 1$).
Corresponding to each monomial $m$, we consider the following
locally-closed subvariety of $\M$:
\begin{equation*}
   \mathfrak Z(m) \defeq \left\{ x\in \M \left|\,
       \lim_{t\to\infty} t\diamond x \in \M(m) \right\}\right.,
\end{equation*}
where $\diamond$ denotes the $\C^*$-action. If $m_P\in\mathcal M'$,
then all monomials appearing in $\qch{q,t}(M(P))$ is contained in
$\mathcal M'$ by \thmref{thm:ind}. Moreover, it can be shown that the
above $\mathfrak Z(m)$ is a lagrangian subvariety of $\M$\footnote{The
class of $\C^*$-actions used for the $\ve$--characters is different
from that for crystals in \cite[\S8]{Na-tensor}. The condition ensures
that the action in the former class is also in the latter class. The
absence of odd cycles is used here.}.
Moving all $m$, the union of all $\mathfrak Z(m)$ forms a closed
lagrangian subvariety $\mathfrak Z$ of $\M$. So the irreducible
components of $\mathfrak Z$ can be identified with monomials. In
\cite{Na-tensor}, a crystal structure is defined on the set of
irreducible components of $\mathfrak Z$, following the work of
Kashiwara-Saito~\cite{KS}.
The crystal structure is isomorphic to a direct sum of the crystals of 
highest weight modules.
We translate this result in terms of monomials, and get the above
theorem.

Applying the above result to $\ve$--characters, we get the following
\begin{Theorem}
  Suppose that $\g$ is of type $ADE$. Let $M(P)$ be a standard module,
  and let $\mathcal M(P)$ be the set of monomials appearing in
  $\qch{q,t}(M(P))$.
Suppose that the monomial $m_P$ corresponding to the 
{\it l\/}--highest weight vector is contained in $\mathcal M'$. 
Then $\mathcal M(P)$ has a structure of a crystal \textup(with respect
to $\g$\textup), which is isomorphic to a direct sum of the crystals
of highest weight modules. Moreover, the crystal graph is obtained from
the graph $\Gamma_P$ by forgetting the multiplicities of monomials and
erasing some arrows.
\end{Theorem}

For example, the crystal of \ref{ex:graph}(3) is
$\B(2\Lambda_1)\oplus\B(\Lambda_2)$ while those of (1), (2) are
crystals of simple modules.

Since the original $\chi_\ve$--character for non-simply-laced case has
a slightly different $A_{i,a}$ from one used in this paper, the proof
(ours or Kashiwara's) of the above theorem does not apply. But the statement
seems to be true.

In general, it is not easy to determine the crystal structure on
$\mathcal M(P)$. But we can do it for a special case.
Choose and fix orientations of edges in the Dynkin diagram. We define
integer $m(i)$ for each vertex $i$ so that $m(i) - m(j) = 1$ if we
have an oriented edge from $i$ to $j$, i.e., $i\to j$. Then we define
$P$ by
\begin{equation*}
   P_i(u) = (1 - ua\ve^{m(i)})^{w_i}
\end{equation*}
for $w_i\in\Z_{\ge 0}$, $a\in\C^*$. A special case is when $M(P)$ is
an {\it l\/}--fundamental representation, i.e., $w_i = 0$ except for
one vertex $i$.

\begin{Proposition}\label{prop:crystal}
Suppose that $\g$ is of type $ADE$ and $P$ as above.
Then the above $\mathcal M(P)$ is isomorphic to the crystal
$\B(\sum_i w_i \Lambda_{i})$ of the highest weight $\g$-module.
\end{Proposition}

This is not true for non-simply-laced case, even if we get the crystal
structure, as conjectured above.

This can be proved by showing that the above lagrangian subvariety
$\mathfrak Z$ is isomorphic to another lagrangian subvariety $\La$,
whose irreducible components are known to be obtained from the highest
weight vector by applying Kashiwara's operators. (The detail depends
on the theory of quiver varieties. So it is not given here.) The first
two examples of \ref{ex:graph} satisfy the assumption of
\propref{prop:crystal}.

\propref{prop:crystal} means that all the monomials appearing in
$\qch{q,t}(M(P))$ can be determined from the crystal $\B(\sum_i w_i
\Lambda_{i})$ of the highest weight $\g$-module. So far, the relation
between the coefficients of monomials and the theory of crystal bases
is unclear. For example, {\it l\/}--fundamental representations are
known to have crystal bases \cite{Kas00}, but their relation to
$\ve$--characters are not known.

\section{type $A_n$}

We number the vertex of the Dynkin graph of type $A_n$ as follows:
\begin{equation*}
\begin{diagram}
  \node{\numbullet{1}}
  \arrow{e,t,-}{}
  \node{\numbullet{2}}
  \arrow{e,t,-}{}
  \node{\cdots}
  \arrow{e,t,-}{}
  \node{\numbullet{n-1}}
  \arrow{e,t,-}{}
  \node{\numbullet{n}}
\end{diagram}
\end{equation*}

We have
\begin{equation*}
   A_{i,a} = Y_{i,a\ve^{-1}} Y_{i,a\ve} Y_{i-1,a}^{-1} Y_{i+1,a}^{-1},
\end{equation*}
where we understand $Y_{i-1,a} = 1$ and $Y_{i+1,a} = 1$ if $i=1$ and
$i=n$ respectively.

We first consider the pullback of the vector representation by the
evaluation homomorphism. It is the {\it l\/}--fundamental
representation $L(\Lambda_1)_a$. Then \propref{prop:crystal} implies
that the vertex of the graph of $\qch{q,t}(L(\Lambda_1)_a)$ is the
same as that of crystal $\B(\Lambda_1)$ and all coefficients are
$1$. Therefore we get
\begin{equation*}
\begin{gathered}[t]
Y_{1,a} \\ || \\ \ffbox{1}_a\!
\end{gathered}
\xrightarrow{1,a\ve}
\begin{gathered}[t]
Y_{1,a\ve^2}^{-1} Y_{2,a\ve} \\ || \\ \ffbox{2}_a\!
\end{gathered}
\xrightarrow{2,a\ve^2}
\dots
\xrightarrow{n\!-\!1,a\ve^{n\!-\!1}}
\begin{gathered}[t]
Y_{n\!-\!1,a\ve^n}^{-1} Y_{n,a\ve^{n\!-\!1}} \\ || \\  \ffbox{n}_a\!
\end{gathered}
\xrightarrow{n,a\ve^n}
\begin{gathered}[t]
Y_{n,a\ve^{n+1}}^{-1} \\ || \\ \ffbox{n\!+\!1}_a\!
\end{gathered}
\end{equation*}
Now it becomes clear that there are no extra arrows. So the graph
$\Gamma_P$ is exactly the same as the crystal graph.

We introduce the symbol $\ffbox{i}_a$ by the above equations. In fact, 
this can be easily shown from \thmref{thm:ind} without any knowledge
about the representation theory.

Let $\mathbf B = \{ 1,\dots,n+1\}$. We give the usual ordering $<$ on
$\mathbf B$.
\begin{Definition}\label{def:tableau}
(1)\ A {\it column tableau\/} $T$ is a map
\begin{equation*}
   T\colon \{ a\ve^{N-1}, a\ve^{N-3}, \dots, a\ve^{1-N} \} \to \mathbf B,
\end{equation*}
for $a\in\C^*$, $1\le N\le n$.
We call $N$ the {\it length\/} and $a$ the {\it center\/} of $T$
respectively.
We associate a monomial $m_T$ to $T$ by
\begin{equation*}
   m_T = \prod_{p=1}^N \ffbox{i_p}_{a\ve^{N+1-2p}},
\end{equation*}
where $i_p = T(a\ve^{N+1-2p})$.
We write this graphically as
\begin{equation*}
T = 
\renewcommand{\arraycolsep}{0pt}
\begin{array}{|c|}
   \hline \hbox to \yl{\hfill$\scriptstyle i_1$\hfill} \\
   \hline \hbox to \yl{\hfill$\vdots$\hfill} \\
   \hline \hbox to \yl{\hfill$\scriptstyle{i_N}$\hfill} \\
   \hline
\end{array}_{a\ve^{1-N}}.
\end{equation*}
For $p\neq N$ the suffix of $i_p$, which is $a\ve^{N+1-2p}$, is
omitted since it can be determined from that of $i_N$.
The same graphical notation will be used for $m_T$.
%
%
We extend $T$ to a map from $\C^*\to \mathbf B\sqcup \{ 0\}$ by
setting
\begin{equation*}
   T(b) = 0 \qquad \text{if $b\neq a\ve^{N-1}, \dots, a\ve^{1-N}$}.
\end{equation*}
In this case, $\{a\ve^{N-1}, \dots, a\ve^{1-N}\}$ will be called the
{\it support\/} of $T$.

(2)\ A {\it tableau\/} $T$ is a finite sequence of column tableaux $T
= (T_1,T_2,\dots,T_L)$.
Its {\it shape\/} is the sequence of lengths and centers of columns:
$(N_1,a_1,N_2,a_2,\dots,N_L,a_L)$.
We write $T$ graphically as
\begin{equation*}
T = 
\renewcommand{\tabcolsep}{0pt}
\begin{tabular}{|c|c|c|c|c|c}
\cline{5-5}
\multicolumn{3}{c}{}&&&
\\
\cline{1-1}
&\multicolumn{2}{c}{}&&&
\\
\cline{2-2}
&&\multicolumn{1}{c}{\dots}&&$T_L$&
\\
&$T_{2}$&\multicolumn{1}{c}{}&&&${}_{a_L\ve^{1-N_L}}$
\\
\cline{5-5}
$T_1$&&\multicolumn{3}{c}{}
\\
\cline{2-2}
&\multicolumn{4}{c}{}&
\\
&\multicolumn{4}{c}{}&
\\
\cline{1-1}
\end{tabular},
\end{equation*}
where $T_\alpha$'s are placed so that $T_1(b)$, $T_2(b)$, \dots,
$T_L(b)$ appear in a row for each $b\in\C^*$.
As above suffixes for $T_\alpha$ ($\alpha\le L-1$) are omitted since they
are determined from the suffix of $T_L$ and the positions of $T_\alpha$.
%
Since we assume that supports of $T_\alpha$ are contained in
$a\ve^{2\Z}$, all rows are matched.

The associated monomial $m_T$ is given by
\begin{equation*}
   m_T = \prod_{\alpha=1}^L m_{T_\alpha}.
\end{equation*}

(3)\ We define $s_{\alpha, \beta}$ by
\begin{equation*}
  \frac{a_\alpha \ve^{N_\alpha}}{a_\beta \ve^{N_\beta}}
  = \ve^{2s_{\alpha,\beta}}.
\end{equation*}
If the left hand side is not in $\ve^{2\Z}$, we simply set
$s_{\alpha,\beta} = - \infty$. This is the number of boxes in
$T_\alpha$ which is located upper than the top of $T_\beta$.
Then we set
\begin{equation*}
\begin{split}
   & d(T_\alpha, T_\beta) \defeq
   \begin{aligned}[t]
   &
   \sum_{b\in\C^*}\left(T_\alpha(b\ve^{-2})<T_\beta(b)<T_\alpha(b)\right)
\\
   & \quad - \left(N_\alpha < T_\beta(a_\alpha q^{-1-N_\alpha}) \le
       T_\alpha(a_\alpha q^{1-N_\alpha})\right)
\\
   & \qquad + \left(N_\alpha - s_{\alpha,\beta} < N_\alpha <
    T_\beta(a_\alpha q^{-1-N_\alpha})\right)
   \end{aligned}
\\
   & d(T) \defeq \sum_{\alpha < \beta} d(T_\alpha,T_\beta).
\end{split}
\end{equation*}
Note that $T_\alpha(a_\alpha q^{1-N_\alpha})$ is the entry of the
bottom of $T_\alpha$ and $T_\beta(a_\alpha q^{-1-N_\alpha})$ is the
entry of $T_\beta$ of one row below.
From the definition, it is clear that $d(T_\alpha, T_\beta) = 0$
unless both of their supports are contained in $a\ve^{2\Z}$ for some
$a\in\C^*$.

(4)\ A tableau $T$ is said to be {\it column increasing\/} if the
entries in each column strictly increase from top to bottom.
\end{Definition}

We have
\begin{equation}\label{eq:m_T}
   m_T = \prod_a \prod_{i=1}^n
      Y_{i,a}^{\# \ffbox{i}_{a\ve^{1-i}}
                    - \# \ffbox{i\!+\!1}_{a\ve^{-1-i}}},
\end{equation}
where $\#\ffbox{i}_a = \#\{ \alpha\mid T_\alpha(a) = i\}$.

\begin{Definition}
Two tableaux $T$ and $T'$ are {\it equivalent\/} if
\(
   \#\{ \alpha\mid T_\alpha(a) = i\}
    =  \#\{ \alpha\mid T'_\alpha(a) = i\}
\)
for all $i=1,\dots, n+1$, $a\in\C^*$.
Namely $T'$ is obtained from $T$ by permuting boxes in the same rows.
\end{Definition}

It is clear that monomials $m_T$ and $m_{T'}$ are equal if $T$ and
$T'$ are equivalent. The converse is not true, but we can determine
when $m_T = m_{T'}$.

\begin{Lemma}\label{lem:mT=mT'}
Let $T$ and $T'$ be tableaux. Then the corresponding monomials $m_T$
and $m_{T'}$ are equal if and only if $T$ and $T'$ become equivalent
after we add several columns of the form
\begin{equation*}
\renewcommand{\arraycolsep}{0pt}
\begin{array}{|c|}
   \hline\hbox to \yl{\hfill$\scriptstyle 1$\hfill}\\
   \hline\hbox to \yl{\hfill$\scriptstyle 2$\hfill}\\
   \hline\hbox to \yl{\hfill$\vdots$\hfill}\\
   \hline\hbox to \yl{\hfill$\scriptstyle{n\!+\!1}$\hfill}\\
   \hline
\end{array}_{a}
\end{equation*}
to $T$ and $T'$.
\end{Lemma}

\begin{proof}
Let
\[
   \#'\ffbox{i}_a \defeq
     \#\{ \alpha\mid T_\alpha(a) = i\}
    -  \#\{ \alpha\mid T'_\alpha(a) = i\}.
\]
By \eqref{eq:m_T} $m_T = m_{T'}$ if and only if
\begin{equation*}
   \#' \ffbox{i}_{a\ve^{1-i}}
                    = \#' \ffbox{i\!+\!1}_{a\ve^{-1-i}}
\end{equation*}
for any $a,i$. Replacing $a$ by $a\ve^{2n+1-i}$, we get
\begin{equation*}
   \#' \ffbox{i}_{a\ve^{2n+2-2i}}
   = \#' \ffbox{i\!+\!1}_{a\ve^{2n-2i}}
\end{equation*}
Namely $m_T$ and $m_{T'}$ are equal if and only if
\begin{equation*}
   \#' \ffbox{i}_{a\ve^{2n+2-2i}}
\end{equation*}
is independent of $i=1,\dots,n+1$ for any $a\in\C^*$.
Let $d_a$ be this integer. If $d_a > 0$, we add $(d_a)$-columns of the
above form to $T'$. If $d_a < 0$, we add $(-d_a)$-columns to $T$. Then 
the resulting tableaux are equivalent.
\end{proof}

\begin{Lemma}\label{lem:ldom}
Let $T$ be a tableau. The corresponding monomial $m_T$ is {\it
l\/}--dominant if and only if $T$ is equivalent to a tableau $T'$ such
that every column is of the form
\begin{equation*}
\renewcommand{\arraycolsep}{0pt}
\begin{array}{|c|}
   \hline \hbox to \yl{\hfill$\scriptstyle 1$\hfill} \\
   \hline \hbox to \yl{\hfill$\scriptstyle 2$\hfill} \\
   \hline \hbox to \yl{\hfill$\vdots$\hfill} \\
   \hline \hbox to \yl{\hfill$\scriptstyle{i}$\hfill} \\
   \hline
\end{array}_{a}
\end{equation*}
for some $a\in\C^*$, $i \in \{ 1, 2,\dots, n+1\}$.
\end{Lemma}

\begin{proof}
For $i=1,\dots n$ we have
\begin{equation*}
   Y_{i,a}
 = \ffbox{1}_{a\ve^{i-1}} \ffbox{2}_{a\ve^{i-3}} \cdots
     \ffbox{i}_{a\ve^{1-i}}
.
\end{equation*}
Thus an {\it l\/}--dominant monomial is given by a tableau $T''$ whose
column is of the form as above with $a\in \C^*$, $i\in \{1,\dots,
n\}$. On the other hand, \lemref{lem:mT=mT'} means that $m_T =
m_{T''}$ if and only if $T$ is equivalent to a tableau $T'$ which is
obtained by adding columns of the above type with $i=n+1$ to $T''$.
\end{proof}

Now we give a tableaux sum expression of $t$--analogs of
$\ve$--characters. We start with {\it l\/}--fundamental representations.

Let
\[
   \B(\Lambda_N)_a = 
   \left.\left\{
   T = 
\renewcommand{\arraycolsep}{0pt}
\begin{array}{|c|}
   \hline \hbox to \yl{\hfill$\scriptstyle i_1$\hfill} \\
   \hline \hbox to \yl{\hfill$\vdots$\hfill} \\
   \hline \hbox to \yl{\hfill$\scriptstyle{i_N}$\hfill} \\
   \hline
\end{array}_{a\ve^{1-N}}
   \right| i_p\in \mathbf B, \; i_1 < i_2 < \dots < i_N \right\}.
\]

\begin{Proposition}\label{prop:Afun}
For $1\le N\le n$, we have
\begin{equation}\label{eq:Afun}
   \qch{q,t}(L(\Lambda_N)_a)
   = \sum_{T\in\B(\Lambda_N)_a} m_T.
\end{equation}
\end{Proposition}

\begin{proof}
We check the conditions~\ref{thm:ind}(1)(2). By \lemref{lem:ldom} it
is clear that the only {\it l\/}--dominant monomial in the right hand
side of \eqref{eq:Afun} is the {\it l\/}--highest weight vector,
i.e., $i_1 = 1$, \dots, $i_N = N$.

Let $T\in \B(\Lambda_N)_a$ as above. The exponent of
$Y_{i,a\ve^{2(i+1-p)}}$ is positive if and only if $i_p = i$, $i_{p+1}
\neq i+1$. In this case the exponent is equal to $1$, and the exponent
of other $Y_{i,b}$'s ($b\neq a\ve^{2(i+1-p)}$) are all $0$. Let $T'$
be the tableau obtained from $T$ by changing $i$ to $i+1$. It is in
$\B(\Lambda_N)_a$ since $i_{p+1} \neq i+1$ and we have
\begin{equation*}
   m_{T} + m_{T'} 
   = Y_{i,a\ve^{2(i+1-p)}}\left(1 + A_{i,a\ve^{2(i-p)+3}}^{-1}\right)M,
\end{equation*}
where $M$ does not contain the factor $Y_{i,b}^\pm$ for any
$b\in\C^*$. This shows that the right hand side of \eqref{eq:Afun}
satisfies the condition~\ref{thm:ind}(1).
\end{proof}

Our next task is to compute
\(
   d(m_{T}, m_P; m_{T'}, m_{P'})
\)
for two column tableaux $T$, $T'$ with corresponding {\it
l\/}--dominant monomials $m_P$, $m_{P'}$.
We represent them graphically as
\begin{equation*}
   T = 
\renewcommand{\arraycolsep}{0pt}
\begin{array}{|c|}
   \hline \hbox to \yl{\hfill$\scriptstyle i_1$\hfill} \\
   \hline \hbox to \yl{\hfill$\vdots$\hfill} \\
   \hline \hbox to \yl{\hfill$\scriptstyle{i_N}$\hfill} \\
   \hline
\end{array}_{a\ve^{1-N}},
\quad
   T' = 
\renewcommand{\arraycolsep}{0pt}
\begin{array}{|c|}
   \hline \hbox to \yl{\hfill$\scriptstyle j_{s\!+\!1}$\hfill} \\
   \hline \hbox to \yl{\hfill$\vdots$\hfill} \\
   \hline \hbox to \yl{\hfill$\scriptstyle{j_{M}}$\hfill} \\
   \hline
\end{array}_{a\ve^{N+1-2M}}.
\end{equation*}
Note that we fix the entries of $T'$ so that $j_p$ and $i_p$ have the
same vertical coordinate if we write $T$ and $T'$ graphically by the
rule \defref{def:tableau}(2).
We set $i_p = j_{p'} = 0$ if $p\neq 1,\dots,N$, $p'\neq s+1,\dots,M$.
The corresponding {\it l\/}--dominant monomials $m_P$ and $m_{P'}$ are
given by $m_P = Y_{N,a}$, $m_{P'} = Y_{M-s,a\ve^{N-M-s}}$.

\begin{Lemma}\label{lem:d}
\begin{equation*}
\begin{split}
   & d(m_{T}, m_P; m_{T'}, m_{P'})\\
=\; & 
  \sum_{p=1}^N (i_{p-1} < j_p < i_p)
  - (N < j_{N+1} \le i_N)
  + (N-s < N < j_{N+1}).
\end{split}
\end{equation*}
\end{Lemma}

\begin{proof}
We have
\begin{equation}\label{eq:v}
\begin{split}
  & m_{T}
  = m_{P} 
     \prod_{p=1}^N \prod_{i=p}^{i_p-1} A_{i,a\ve^{N+1-2p+i}}^{-1},
\\
  & m_{T'}
  = m_{P'} 
     \prod_{p=s+1}^M \prod_{i=p-s}^{j_p-1} A_{i,a\ve^{N+1-2p+i}}^{-1}.
\end{split}
\end{equation}
Hence we have
\begin{equation*}
   v_{i,a\ve^{N+1-2p+i}}(m_T,m_P)
   = (p\le i \le i_p-1).
\end{equation*}
On the other hand, we have
\begin{equation*}
   u_{i,a\ve^{N-2p+i}}(m_{T'})
   = (j_p = i) - (j_{p+1}=i+1).
\end{equation*}
Thus we get
\begin{equation*}
\begin{split}
   & \sum_i v_{i,a\ve^{N-2p+i+1}}(m_T,m_P) u_{i,a\ve^{N-2p+i}}(m_{T'})
\\
   = \; &(p \le j_p \le i_p-1) - (p+1\le j_{p+1} \le i_p).
\end{split}
\end{equation*}
Summing up with respect to $p$, we get
\begin{equation*}
\begin{split}
   & \sum_{p=1}^N
   \sum_i v_{i,a\ve^{N-2p+i+1}}(m_T,m_P) u_{i,a\ve^{N-2p+i}}(m_{T'})
\\
   = \; &
   \sum_{p=1}^N (p \le j_p \le i_p-1) - (p+1\le j_{p+1} \le i_p)
\\
   = \; &
   \sum_{p=1}^{N} \left[(p \le j_p \le i_p-1)
      - (p\le j_{p} \le i_{p-1})\right] - (N+1\le j_{N+1} \le i_N)
\end{split}
\end{equation*}
Here we have used that $p\le j_{p} \le i_{p-1}$ never hold if $p=1$ in
the last equality.
Note that $i_{p-1} < i_p$ for $p=1,\dots, N$. Thus
\(
   p\le j_{p} \le i_{p-1}
\)
implies
\(
  p \le j_p \le i_p-1
\).
Hence each term of the above summation is
\begin{equation*}
   (i_{p-1} < j_p < i_p).
\end{equation*}
Note that $p \le j_p$ holds automatically since $p-1 \le i_{p-1}$.

We have $u_{N,a}(m_{P^1}) = 1$ and other $u_{i,b}(m_{P^1})$'s are all
$0$. By \eqref{eq:v} we have
\begin{equation*}
   v_{N,a\ve^{-1}}(m_{T'},m_{P'}) = (N+1-s \le N < j_{N+1}).
\end{equation*}
Combining all together, we get the assertion.
\end{proof}

Now let $M(P)$ be an arbitrary standard module. We decompose $P = P^1
P^2 P^3\cdots$ so that each $M(P^\alpha)$ is an {\it l\/}--fundamental
representation and the condition~\ref{eq:Z} is met with respect to the
ordering.
(There might be several orderings satisfying \eqref{eq:Z}. In that
case, we just fix one such ordering.)
Let $(N_\alpha, a_\alpha)$ be the shape of a column tableau
corresponding to $M(P^\alpha)$ by \propref{prop:Afun}.
Then let $\B(P)$ be the set of column increasing tableaux with shape
$(N_1, a_1, N_2, a_2, \dots, N_L, a_L)$.
We apply \thmref{thm:ind}(3) successively to get the following:
\begin{Theorem}
\begin{equation*}
    \qch{q,t}(M(P))
    = \sum_{T\in\B(P)} t^{2d(T)} m_T.
\end{equation*}
\end{Theorem}

\begin{Example}
Let $\g$ be of type $A_2$. We give only Young tableaux. The
corresponding $m_T$ and $d(T)$ are given in Example~\ref{ex:graph}.

(1) $\qch{q,t}(M(P))$ with $m_P = \Y 11 \Y 2\ve$ is given by
\newcommand{\OneTwoY}[3]{
\renewcommand{\arraycolsep}{0pt}
\begin{array}{cc}
   \cline{2-2}
   & \vline\hbox to \yl{\hfill$\scriptstyle #1$\hfill}\vline
\\
   \cline{1-2}
   \vline\hbox to \yl{\hfill$\scriptstyle #2$\hfill} &
   \vline\hbox to \yl{\hfill$\scriptstyle #3$\hfill}\vline
\\
   \cline{1-2}
\end{array}_1
}
\begin{equation*}
\begin{CD}
  \OneTwoY 112
  @>{2,\ve^2}>>
  \OneTwoY 113
  @>{1,\ve^3}>>
  \OneTwoY 213
  @>{1,\ve}>>
  \OneTwoY 223
\\
  @V{1,\ve}VV @V{1,\ve}VV @. @VV{2,\ve^2}V
\\
  \OneTwoY 122
  @>>{2,\ve^2}>
  \OneTwoY 123 + \OneTwoY 132
  @>>{2,\ve^2}>
  \OneTwoY 133
  @>>{1,\ve^3}>
  \OneTwoY 233
\end{CD}
\end{equation*}

(2) $\qch{q,t}(M(P))$ with $m_P = \Y 11^2$ is given by
\newcommand{\TwoY}[2]{
\renewcommand{\arraycolsep}{0pt}
\begin{array}{cc}
   \cline{1-2}
   \vline\hbox to \yl{\hfill$\scriptstyle #1$\hfill} &
   \vline\hbox to \yl{\hfill$\scriptstyle #2$\hfill}\vline
\\
   \cline{1-2}
\end{array}_1
}
\begin{equation*}
  \begin{CD}
  \TwoY 11
  @>{1,\ve}>>
  \TwoY 12 + \TwoY 21
  @>{1,\ve}>>
  \TwoY 22
  @.
\\
  @. @V{2,\ve^2}VV @V{2,\ve^2}VV
\\
  @. 
  \TwoY 13 + \TwoY 31
   @>>{1,\ve}>
  \TwoY 23 + \TwoY 32
  @>>{2,\ve^2}>
  \TwoY 33
\end{CD}
\end{equation*}

(3) $\qch{q,t}(M(P))$ with $m_P = \Y 11 \Y 1{\ve^2}$ is given by
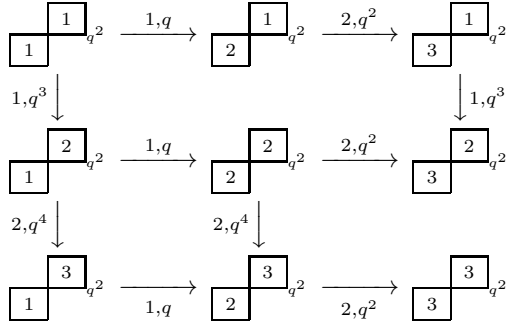
\newcommand{\OneOneY}[2]{
\renewcommand{\arraycolsep}{0pt}
\begin{array}{c|c}
   \cline{2-2}
   & \hbox to \yl{\hfill$\scriptstyle #1$\hfill}\vline
\\
   \cline{1-2}
   \vline\hbox to \yl{\hfill$\scriptstyle #2$\hfill} &
\\
   \cline{1-1}
\end{array}{\scriptscriptstyle \ve^2}
}
\begin{equation*}
\begin{CD}
  \OneOneY 11 @>{1,\ve}>> \OneOneY 12 @>{2,\ve^2}>> \OneOneY 13
\\
  @V{1,\ve^3}VV @. @VV{1,\ve^3}V
\\
  \OneOneY 21 @>{1,\ve}>> \OneOneY 22 @>{2,\ve^2}>> \OneOneY 23
\\
  @V{2,\ve^4}VV @V{2,\ve^4}VV @.
\\
  \OneOneY 31 @>>{1,\ve}> \OneOneY 32 @>>{2,\ve^2}> \OneOneY 33
\end{CD}
\end{equation*}
\end{Example}

\begin{Remark}
(1) In \cite{NT} a different convention for tableaux was used. Each
column is located so that $T_1(b)$, $T_2(b\ve^2)$, \dots appear in the
same row.

(2) When the shapes of tableaux are those of ordinary Young tableaux,
i.e., the tops of columns are the same and the lengths are
nonincreasing, $d(T)$ is equal to $C - l(T^t)$ where
$T^t$ is the transpose of $T$,
$l(T^t)$ was defined in \cite{Na:Hom},
and $C$ is a constant depending only the shape of $T$ and numbers of
figures.
(In fact, $C$ is equal to $\frac12 \dim_\C \M(\bv,\bw)$, where
$\M(\bv,\bw)$ is the quiver variety containing the point corresponding
to $T$.)
Moreover, the assumption of \propref{prop:crystal} is satisfied in
this case. (See also \cite[\S8.5]{Na-qchar}.)
\end{Remark}

\section{type $D_n$}

We number the vertices as follows.
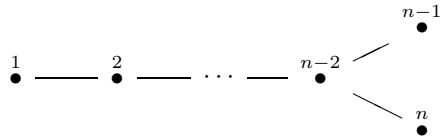
\begin{equation*}
\dgARROWLENGTH=1.5em
\begin{diagram}
  \node[10]{\numbullet{n-1}}
\\
  \node[2]{\numbullet{1}}
  \arrow[2]{e,t,-}{}
  \node[2]{\numbullet{2}}
  \arrow[2]{e,t,-}{}
  \node[2]{\cdots}
  \arrow[2]{e,t,-}{}
  \node[2]{\numbullet{n-2}}
  \arrow{ene,t,-}{}
  \arrow{ese,t,-}{}
\\
  \node[10]{\numbullet{n}}
\end{diagram}
\end{equation*}

The vector representation of $\g$ is known to be lifted to a
$\Ul$-module. It is an {\it l\/}--fundamental representation
$L(\Lambda_1)_a$. The graph of its $\ve$-character is the same as that 
of crystal $\B(\Lambda_1)$ as in the case of $A_n$. It is
\begin{equation*}
\dgARROWLENGTH=1.5em
\dgHORIZPAD=1em
\dgVERTPAD=1ex
\renewcommand{\dgeverylabel}{\scriptscriptstyle}
\begin{diagram}
\node[5]{\ffbox{n}_a\!}
\arrow{se,t}{n,a\ve^{n\!-\!1}}
\\
\node{\ffbox{1}_a\!}
\arrow{e,t}{1,a\ve}
\node{\ffbox{2}_a\!}
\arrow{e,t}{2,a\ve^2}
\node{\cdots}
\arrow{e,t}{\!\!n\!-\!2,a\ve^{n\!-\!2}}
\node{\ffbox{n\!-\!1}_a\!}
\arrow{ne,l}{n\!-\!1,a\ve^{n\!-\!1}}
\arrow{se,b}{n,a\ve^{n\!-\!1}}
\node[2]{\ffbox{\overline{n\!-\!1}}_a\!}
\arrow{e,t}{n\!-\!2,a\ve^n}
\node{\cdots}
\arrow{e,t}{2,a\ve^{2n\!-\!4}}
\node{\ffbox{\overline{2}}_a\!}
\arrow{e,t}{\!1,a\ve^{2n\!-\!3}}
\node{\ffbox{\overline{1}}_a\!}
\\
\node[5]{\ffbox{\overline{n}}_a\!}
\arrow{ne,b}{n\!-\!1,a\ve^{n\!-\!1}}
\end{diagram}
,
\end{equation*}
where
\begin{equation*}
\begin{aligned}[c]
& \ffbox{i}_a = Y_{i-1,a\ve^i}^{-1} Y_{i,a\ve^{i-1}} \qquad(1 \le i \le n-2)
\\
& \ffbox{n\!-\!1}_a
 = Y_{n-2,a\ve^{n-1}}^{-1} Y_{n-1,a\ve^{n-2}} Y_{n,a\ve^{n-2}} 
\\
& \ffbox{n}_a = Y_{n-1,a\ve^{n}}^{-1} Y_{n,a\ve^{n-2}}
\\
& \ffbox{\overline{n}}_a = Y_{n-1,a\ve^{n-2}} Y_{n,a\ve^{n}}^{-1}
\\
& \ffbox{\overline{n\!-\!1}}_a = Y_{n-2,a\ve^{n-1}} Y_{n-1,a\ve^{n}}^{-1}
Y_{n,a\ve^{n}}^{-1}
\\
& \ffbox{\overline{i}}_a = Y_{i-1,a\ve^{2n-2-i}} Y_{i,a\ve^{2n-1-i}}^{-1}
\qquad(1 \le i \le n-2).
\end{aligned}
\end{equation*}
Here $Y_{0,b}$ is understood as $1$. This is also easily shown by
\thmref{thm:ind}. (The notation is borrowed from \cite{Kas-Na}.)

Let $\mathbf B = \{ 1,\dots,n,\overline{n},\dots, \overline{1}\}$. 
We give the ordering $\prec$ on the set $\mathbf B$
by
\begin{equation*}
  1 \prec 2 \prec \cdots \prec n-1 \prec
  \begin{matrix}n \\ \\ \overline{n}\end{matrix}
 \prec \overline{n-1} \prec\cdots \prec \overline{2}\prec \overline{1}.
\end{equation*}
Remark that there is no order between $n$ and $\overline{n}$.

We define a tableau $T$ and its associated monomial exactly as in the
type $A_n$ case. We just replace $\mathbf B$. (So far we do not
include column corresponding to spin representations.)
The following is an analog of \lemref{lem:mT=mT'}. (In fact, it will
not be used later.)

\begin{Lemma}\label{lem:D:mT=mT'}
Let $T$ and $T'$ be tableaux. Then the corresponding monomials $m_T$
and $m_{T'}$ are equal if and only if $T$ and $T'$ become equivalent
after we add several pairs of columns
\begin{equation*}
   \begin{array}{|c|}
   \hline
   \scriptstyle 1 \\
   \hline \vdots \\
   \hline \scriptstyle{i} \\ \hline
   \end{array}_{a}
,\quad
   \begin{array}{|c|}
   \hline
   \scriptstyle \overline{i} \\ \hline \vdots \\
   \hline \scriptstyle \overline{1} \\ \hline
   \end{array}_{a\ve^{2-2n}}
\end{equation*}
for some $i\in \{1,\dots, n\}$, $a\in\C^*$ to $T$ and $T'$.
\end{Lemma}

\begin{proof}
Let $\#'\ffbox{i}_a$ be the number of boxes with entry $i$ in the
row corresponding to $a$ of $T$ minus that of $T'$.
Then $m_T = m_{T'}$ if and only if
\begin{equation*}
\begin{split}
    & \#'\ffbox{i}_{a\ve^{1-i}} - \#'\ffbox{i+1}_{a\ve^{-1-i}}
    = \#'\ffbox{\overline{i}}_{a\ve^{i-2n+1}}
    -  \#'\ffbox{\overline{i\!+\!1}}_{a\ve^{i-2n+3}} \quad
    \text{for $1\le i\le n-2$} \\
    & \#'\ffbox{n\!-\!1}_{a\ve^{2-n}} - \#'\ffbox{n}_{a\ve^{-n}}
    = \#'\ffbox{\overline{n\!-\!1}}_{a\ve^{-n}}
    -  \#'\ffbox{\overline{n}}_{a\ve^{2-n}}
\\
    & \#'\ffbox{n\!-\!1}_{a\ve^{2-n}} + \#'\ffbox{n}_{a\ve^{2-n}}
    = \#'\ffbox{\overline{n\!-\!1}}_{a\ve^{-n}}
    + \#'\ffbox{\overline{n}}_{a\ve^{-n}}.
\end{split}
\end{equation*}
From the second and third equations we get
\begin{equation*}
   \#'\ffbox{n}_{a\ve^{2-n}} 
     + \#'\ffbox{n}_{a\ve^{-n}}
   = \#'\ffbox{\overline{n}}_{a\ve^{2-n}}
     + \#'\ffbox{\overline{n}}_{a\ve^{-n}}.
\end{equation*}
Moving $a\in\C^*$, we get
\begin{equation*}
   \#'\ffbox{n}_{a}
   = \#'\ffbox{\overline{n}}_{a}.
\end{equation*}
Set this number $d_{n,a}$. Substituting this back to the second
equation, we get
\begin{equation*}
   \#'\ffbox{n\!-\!1}_{a\ve^2} - d_{n,a}
    = \#'\ffbox{\overline{n\!-\!1}}_{a}
    - d_{n,a\ve^2}.
\end{equation*}
Set this number $d_{n-1,a\ve^2}$. Using the first equation, we define
$d_{i,a\ve^{2(n-i)}}$ inductively by
\begin{equation*}
\begin{split}
   & d_{i,a\ve^{2(n-i)}} 
\\
   \defeq \;& \#' \ffbox{i}_{a\ve^{2(n-i)}}
     - \sum_{j: i+1\le j\le n} d_{j,a\ve^{2(n-j)}}
   = \#' \ffbox{\overline{i}}_{a} 
          - \sum_{j:i+1\le j\le n} d_{j,a\ve^{2(n-i)}}.
\end{split}
\end{equation*}
For each $i\in \{ 1,\dots, n\}$, we add $(d_{i,a})$-pairs of columns
(as in the statement) to $T'$ if $d_{i,a} > 0$ and we add
$(-d_{i,a})$-pairs to $T$ if $d_{i,a} < 0$.
The resulting tableaux are equivalent.
\end{proof}

We also have an obvious analog of \lemref{lem:ldom} for $D_n$.

Let
\[
   \B(\Lambda_N)_a = 
   \left\{\left.
   T = 
\renewcommand{\arraycolsep}{0pt}
\begin{array}{|c|}
   \hline \hbox to \yl{\hfill$\scriptstyle i_1$\hfill} \\
   \hline \hbox to \yl{\hfill$\vdots$\hfill} \\
   \hline \hbox to \yl{\hfill$\scriptstyle{i_N}$\hfill} \\
   \hline
\end{array}_{a\ve^{1-N}}
   \right| i_p\in \mathbf B,
   \; i_1 \nsucceq i_2 \nsucceq \dots \nsucceq i_N \right\}.
\]
We define the associated degree by
\begin{equation}\label{eq:expo}
  l(T) = \# \{ p \mid \text{$i_p = i, i_{p+n-1-i} = \overline{i}$ for
  $i=1,\dots, n-2$}\}.
\end{equation}

For $i = 1, \dots, n$ and an integer $s$, we define $p(i,s)$ and
$p'(i,s)$ so that
\begin{equation*}
\begin{cases}
   s = N - 2p(i,s) + i = N - 2p'(i,s) - 2 + 2n - i
   & \text{if $1\le i \le n-1$},
\\
   s = N - 2p(n,s) + n-1 = N - 2p'(n,s) + n + 1
   & \text{if $i = n$}.
\end{cases}
\end{equation*}
If such $p(i,s)$, $p'(i,s)$ do not exist, they are undefined. We have
\begin{equation}\label{eq:ii}
   u_{i,a\ve^s}(m_T)
  = 
  \begin{cases}
  \begin{aligned}[c]
  & \left(i_{p(i,s)} = i\right) - \left(i_{p(i,s)+1} = i+1\right)
  \\
  &\qquad + \left(i_{p'(i,s)} = \overline{i+1}\right)
      - \left(i_{p'(i,s)+1} = \overline{i}\right)
  \end{aligned}
   & \text{if $1\le i\le n-1$},
\\
  \begin{aligned}[c]
    & \left(i_{p(n,s)} = n-1\right) + \left(i_{p(n,s)} = n\right)
    \\
    & \qquad
      - \left(i_{p'(n,s)} = \overline{n}\right)
      - \left(i_{p'(n,s)} = \overline{n-1}\right)
  \end{aligned}
   & \text{if $i = n$}.
  \end{cases}
\end{equation}
If $p(i,s)$, $p'(i,s)$ are undefined, the right hand side is
understood as $0$.

We also have
\begin{equation}\label{eq:v1}
   v_{i,aq^{s+1}}(m_T, Y_{N,a}) =
   \begin{cases}
   \left(p(i,s) \le i \prec i_{p(i,s)}\right)
      + \left(\overline{i} \preceq i_{p'(i,s)}\right)
      & \text{if $1\le i\le n-2$},
   \\
   \left(n \preceq i_{p(n-1,s)}\right) & \text{if $i=n-1$},
   \\
   \left(\overline{n} \preceq i_{p(n,s)}\right)
   & \text{if $i=n$}.
   \end{cases}
\end{equation}

\begin{Proposition}\label{prop:Dfun}
For $1\le N\le n-2$, we have
\begin{equation}\label{eq:qchar}
    \qch{q,t}(L(\Lambda_N)_a)
    = \sum_{T\in \B(\Lambda_N)_a} t^{2l(T)} m_T.
\end{equation}
\end{Proposition}

\begin{proof}
Let $m$ be a monomial in $Y_{i,a\ve^s}^\pm$ for fixed $a$.
Following \cite{FM}, we say $m$ is {\it right negative\/} if the
factor $Y_{i,a\ve^s}$ appearing in $m$, for which $s$ is maximal, have
negative powers. The product of right negative monomials is right
negative. An {\it l\/}-dominant monomial is not right negative.

Let us prove that if $m_T$ is {\it not\/} right negative,
then $i_1 = 1$, \dots, $i_N = N$ by induction.
Since $m_T$ is {\it not\/} right negative, there exists $p_0$ such that
$\ffbox{i_{p_0}}_{a\ve^{N+1-2p_0}}$ is {\it not\/} right negative,
that is $i_{p_0} = 1$.
By the rule $i_p \nsucceq i_{p+1}$, we have $p_0 = 1$, i.e.,
$i_1 = 1$. This proves the first step of the induction.

Suppose that we know $i_p = p$ for $p = 1, \dots, k$.
Consider
\[
   m' = m_T Y_{k,a\ve^{N-k}}^{-1}
  = \ffbox{i_{k+1}}_{a\ve^{N-1-2k}}
  \ffbox{i_{k+2}}_{a\ve^{N-3-2k}}\cdots 
  \ffbox{i_N}_{a\ve^{1-N}}.
\]
This is right negative since all $\ffbox{i_p}_{a\ve^{N+1-2p}}$
($p=k+1,\dots, N$) are so. Therefore, the factor $Y_{i,a\ve^s}$
appearing in $m'$, for which $s$ is maximal, must be equal to
$Y_{k,a\ve^{N-k}}^{-1}$ since $m_T$ is not right negative.
By \eqref{eq:ii}, $Y_{k,a\ve^{N-k}}^{-1}$ appears only when
$i_{p(k,N-k)+1} = k+1$ or $i_{p'(k,N-k)+1} = \overline{k}$. We have
\(
   p(k,N-k) = k,
\)
\(
   p'(k,N-k) = n-1.
\)
So the latter case does not occur since $p'(k,N-k)\le N \le n-2$. Thus we
have $i_p = p$ with $p = k+1$.
This completes the induction step. In particular, the only {\it
l\/}--dominant term in \eqref{eq:qchar} is the {\it l\/}--highest weight
term $Y_{N,a}$.

Next we show that the right hand side of \eqref{eq:qchar} satisfies
the condition \ref{thm:ind}(1) for $i=1,\dots,n$. We only give the
proof for the case $i\neq n$. The case $i = n$ can be checked in a
similar way.

Let $T$ be as above. We consider the following statement for $T$:
\begin{enumerate}
\item $i$ occurs, but $i+1$ does not occur.
\item $i$ does not occur, but $i+1$ occur.
\item Both $i$ and $i+1$ occur (hence consecutively), or neither
occurs.
\item $\overline{i+1}$ occurs, but $\overline{i}$ does not occur.
\item $\overline{i+1}$ does not occur, but $\overline{i}$ occur.
\item Both $\overline{i+1}$ and $\overline{i}$ occur (hence
consecutively), or neither occurs.
\end{enumerate}
Tableaux of type (1) and (2) appear in pairs, i.e., they are obtained
by the replacement of $i$ and $i+1$. If $T$ and $T'$ are such a pair,
we have
\begin{equation}\label{eq:match}
   m_T + m_{T'} = Y_{i,a\ve^s} (1 + A_{i,a\ve^{s+1}}^{-1})M,
\end{equation}
for some $s$. Here $M$ is the contribution from the other terms.
Similarly sequences of type (4) and (5) appear in pairs, and we have
\eqref{eq:match} for the pair $(T,T')$. Monomials of sequences of other
types (i.e., (3) and (6)) do not contain $Y_{i,a\ve^s}^\pm$. This
proves the assertion when $t=1$.

Our remaining task is to study the exponent of $t$.
First consider the case when $T$ satisfies both (1) and (4). So
\(
   i_p = i, i_{p'} = \overline{i+1}
\)
for some $p$, $p'$.
Then $T$ is a member of a quadruplet $(T,T',T'',T''')$, where other
members are obtained from $T$ by replacing $i$, $\overline{i+1}$ by
$i+1$, $\overline{i}$.
First consider the case $p' = p+n-1-i$. We have
\begin{equation*}
\begin{split}
   & t^{2l(T)}m_{T} + t^{2l(T')}m_{T'}
      + t^{2l(T'')}m_{T''} + t^{2l(T''')}m_{T'''}
\\
=\; & Y_{i,a\ve^{N-2p+i}}^2
   \left(1 + (t^2 + 1)A_{i,a\ve^{N-2p+i+1}}^{-1} 
     + A_{i,a\ve^{N-2p+i+1}}^{-2}\right)M
\\
=\; & E_i(Y_{i,a\ve^{N-2p+i}}^2  M)
\end{split}
\end{equation*}
by the definition of $l(T)$ in \eqref{eq:expo}. Here $M$ does not contain
the factor $Y_{i,b}^\pm$ for any $b\in\C^*$. This contribution
satisfies the condition \ref{thm:ind}(1).

Next consider the case $p' = p+n-2-i$. We have
\begin{equation*}
\begin{split}
   & t^{2l(T)}m_{T} + t^{2l(T')}m_{T'}
      + t^{2l(T'')}m_{T''} + t^{2l(T''')}m_{T'''}
\\
=\; & Y_{i,a\ve^{N-2p+k}} Y_{i,a\ve^{N-2p+2+k}}
\\
    & \qquad
   \times \left(1 + t^2 A_{i,a\ve^{N-2p+i+1}}^{-1}
     + A_{i,a\ve^{N-2p+i+3}}^{-1}
     + A_{i,a\ve^{N-2p+i+1}}^{-1} A_{i,a\ve^{N-2p+i+3}}^{-1}\right)M
\\
=\; &  Y_{i,a\ve^{N-2p+i}} Y_{i,a\ve^{N-2p-2+i}}
   \left(1 + A_{i,a\ve^{N-2p+i+1}}^{-1}\right)
   \left(1 + A_{i,a\ve^{N-2p+i-1}}^{-1}\right)M
\\
   & \qquad + (t^2-1)Y_{i,a\ve^{N-2p+i}}
    Y_{i,a\ve^{N-2p-2+i}}A_{i,a\ve^{N-2p+i+1}}^{-1}M.
\end{split}
\end{equation*}
Since 
\[
   Y_{i,a\ve^{N-2p+i}} Y_{i,a\ve^{N-2p+2+i}}A_{i,a\ve^{N-2p+i+1}}^{-1}
   \in \Z[Y_{j,b}^\pm]_{\substack{j: j\neq i\\ b\in\C^*}},
\]
this contribution also satisifies \ref{thm:ind}(1).

In the remaining case $p'\neq p+n-1-i$, $p+n-2-i$, we have
\begin{equation*}
\begin{split}
  & t^{2l(T)}m_{T} + t^{2l(T')}m_{T'}
      + t^{2l(T'')}m_{T''} + t^{2l(T''')}m_{T'''}
\\
  =\; & Y_{i,a\ve^{N-2p+i}} Y_{i,a\ve^{N-2p'-2+2n-i}}
  \left(1+A_{i,a\ve^{N-2p+i+1}}^{-1}\right)
  \left(1+A_{i,a\ve^{N-2p'-1+2n-i}}^{-1}\right)M.
\end{split}
\end{equation*}
This also satisfies \ref{thm:ind}(1).

Next consider the case when $T$ satisfies both (1) and (6). We have
$i_p = i$ for some $p$. Then $T$ appears in a pair $(T,T')$ as
above. The exponents $l(T)$ and $l(T')$ are possibly different only if
$i_{p+n-1-i} = \overline{i}$ or $i_{p+n-2-i} = \overline{i+1}$. But
the condition (6) implies $i_{p+n-2-i} = \overline{i+1}$ {\it and\/}
$i_{p+n-1-i} = \overline{i}$ in either cases. Thus the exponents
$l(T)$ and $l(T')$ are the same even in this case. Then we have
\begin{equation*}
   t^{2l(T)}m_{T} + t^{2l(T')}m_{T'}
   = Y_{i,a\ve^{N-2p+i}} (1 + A_{i,a\ve^{N-2p+i+1}}^{-1})M.
\end{equation*}
This satisfies \ref{thm:ind}(1). In the case $T$ satisifies (3) and
(4), we similarly have \ref{thm:ind}(1). In the remaining case when
$T$ satisfies (3) and (6), $m_{T}$ does not contain the factor
$Y_{i,b}$ for any $b$. Thus it satisfies \ref{thm:ind}(1). Hence the
right hand side of \eqref{eq:qchar} satisfies \ref{thm:ind}(1).
\end{proof}

\begin{Example}\label{exm:D4}
Let $\g = D_4$ and $M(P) = L(\Lambda_2)_{1}$. The graph $\Gamma_P$ is
Figure~\ref{fig:D4}. The same example was appeared in
\cite[5.3.2]{Na-qchar}. The subscripts $\ve^{-1}$ are not written.
\begin{figure}[p]
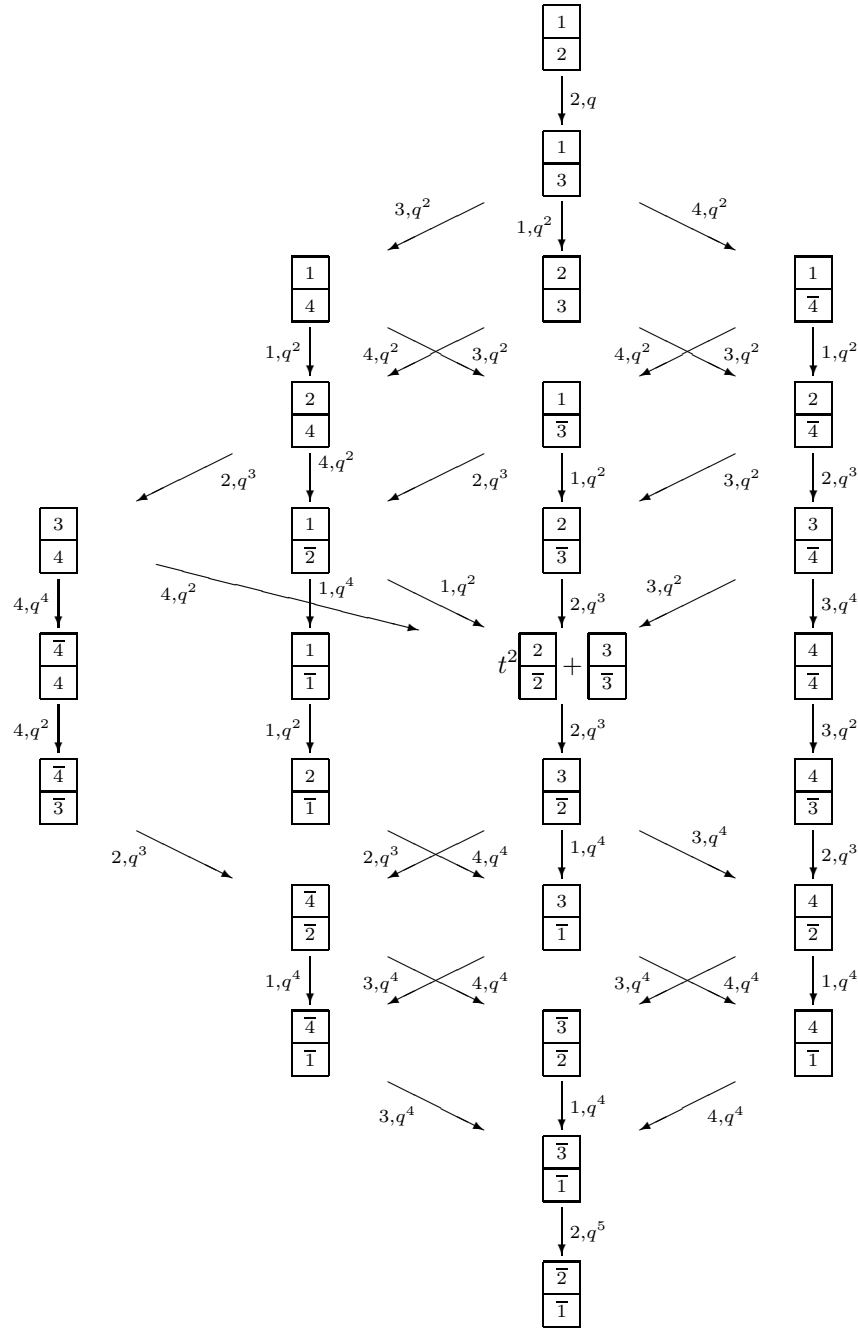

\begin{center}
\leavevmode
\begin{equation*}
\dgARROWLENGTH=0.8em
\dgARROWPARTS=6
\dgHORIZPAD=6em
\dgVERTPAD=1ex
\begin{diagram}
  \node[3]{\twoY 1 2} \arrow{s,r}{2,\ve}
\\
  \node[3]{\twoY 1 3}
  \arrow{sw,t,3}{3,\ve^2} \arrow{s,l,3}{1,\ve^2}
  \arrow{se,t,3}{4,\ve^2}
\\
  \node[2]{\twoY 1 4}
  \arrow{s,l}{1,\ve^2} \arrow{se,b,1}{4,\ve^2}
  \node{\twoY 2 3}
  \arrow{sw,b,1}{3,\ve^2}\arrow{se,b,1}{4,\ve^2}
  \node{\twoY 1 {\overline{4}}}
  \arrow{sw,b,1}{3,\ve^2} \arrow{s,r}{1,\ve^2}
\\
  \node[2]{\twoY 2 4}
  \arrow{sw,b,1}{2,\ve^3} \arrow{s,b,1}{4,\ve^2}
  \node{\twoY 1 {\overline{3}}}
  \arrow{sw,b,1}{2,\ve^3}
  \arrow{s,r}{1,\ve^2}
  \node{\twoY 2 {\overline{4}}}
  \arrow{sw,b,1}{3,\ve^2}
  \arrow{s,r}{2,\ve^3}
\\
  \node{\twoY 3 4}
  \arrow{s,l}{4,\ve^4}
  \arrow{ese,b,1}{4,\ve^2}
  \node{\twoY 1 {\overline{2}}}
  \arrow{s,r,1}{1,\ve^4}
  \arrow{se,t,3}{1,\ve^2}
  \node{\twoY 2{\overline{3}}}
  \arrow{s,r}{2,\ve^3}
  \node{\twoY 3{\overline{4}}}
  \arrow{sw,t,3}{3,\ve^2}
  \arrow{s,r}{3,\ve^4}
\\
  \node{\twoY {\overline{4}}4}
  \arrow{s,l}{4,\ve^2}
  \node{\twoY 1{\overline 1}}
  \arrow{s,l}{1,\ve^2}
  \node{t^2\twoY 2{\overline 2} + \twoY 3{\overline 3}}
  \arrow{s,r}{2,\ve^3}
  \node{\twoY 4{\overline 4}}
  \arrow{s,r}{3,\ve^2}
\\
  \node{\twoY {\overline 4}{\overline 3}}
  \arrow{se,b,1}{2,\ve^3}
  \node{\twoY 2{\overline 1}}
  \arrow{se,b,1}{2,\ve^3}
  \node{\twoY 3{\overline 2}}
  \arrow{sw,b,1}{4,\ve^4}\arrow{s,r,2}{1,\ve^4}
  \arrow{se,t,3}{3,\ve^4}
  \node{\twoY 4{\overline 3}}
  \arrow{s,r}{2,\ve^3}
\\
  \node[2]{\twoY {\overline 4}{\overline 2}}
  \arrow{s,l}{1,\ve^4}
  \arrow{se,b,1}{3,\ve^4}
  \node{\twoY 3{\overline 1}}
  \arrow{sw,b,1}{4,\ve^4}\arrow{se,b,1}{3,\ve^4}
  \node{\twoY 4{\overline 2}}
  \arrow{sw,b,1}{4,\ve^4} \arrow{s,r}{1,\ve^4}
\\
  \node[2]{\twoY {\overline 4}{\overline 1}}
  \arrow{se,b,2}{3,\ve^4}
  \node{\twoY {\overline 3}{\overline 2}}
  \arrow{s,r}{1,\ve^4}
  \node{\twoY 4{\overline 1}}
  \arrow{sw,b,2}{4,\ve^4}
\\
  \node[3]{\twoY {\overline 3}{\overline 1}}
  \arrow{s,r}{2,\ve^5}
\\
  \node[3]{\twoY {\overline 2}{\overline 1}}
\end{diagram}
\end{equation*}
\caption{The graph for $L(\Lambda_2)_{1}$}
\label{fig:D4}
\end{center}
\end{figure}
A careful reader finds that the tableaux appearing here are
slightly different from those in \cite{Kas-Na}.
\end{Example}

\begin{Remark}
Let $\operatorname{Res}$ be the functor sending $\Ul$-modules to
$\Uq$-modules by restriction.
As was shown in \cite{KuS}, \propref{prop:Dfun} implies that the
restriction $\operatorname{Res}L(\Lambda_N)_a$ decomposes as
\[
  \operatorname{Res}L(\Lambda_N)_a =
  \begin{cases}
  V(\Lambda_N)\oplus V(\Lambda_{N-2})\oplus \cdots\oplus
  V(\Lambda_3)\oplus V(\Lambda_1)
  & \text{if $N$ is odd},
  \\
  V(\Lambda_N)\oplus V(\Lambda_{N-2})\oplus \cdots\oplus
  V(\Lambda_2)\oplus V(0)
  &  \text{if $N$ is even},
  \end{cases}
\]
where $V(\lambda)$ is the irreducible highest weight $\Uq$-module with 
highest weight $\lambda$.
This follows from the observation that the (ordinary) character of
$V(\Lambda_N)$ is also described by tableaux sum, but with an extra
condition that $(i_p, i_{p+1})\neq (\overline{n}, n)$.
\end{Remark}

\subsection{Spin representations}

It is known that spin representations of $\g$ can be lifted to a
$\Ul$-module. As in above cases, the vertex of the graph of
$\qch{q,t}$ is the same as that of crystal and all coefficients are
$1$. Following \cite{Kas-Na}, we introduce the half size numbered box:
\begin{equation*}
\begin{split}
   \fhbox{i}_a & = 
   \begin{cases}
     Y_{i-1,a\ve^{i-1}}^{-1} Y_{i,a\ve^{i-2}} 
       & \text{if $1\le i\le n-2$},
   \\
     Y_{n-2,a\ve^{n-2}}^{-1} & \text{if $i=n-1$},
   \\
     Y_{n,a\ve^{n-1}} & \text{if $i=n$},
   \end{cases}
\\
   \fhbox{\overline{i}}_a &= 
   \begin{cases}
     1 & \text{if $1\le i\le n-2$},
   \\
   Y_{n-1,a\ve^{n+1}}^{-1}Y_{n,a\ve^{n+1}}^{-1} 
     & \text{if $i=n-1$},
   \\
   Y_{n-1,a\ve^{n-1}} & \text{if $i=n$}.
   \end{cases}
\end{split}
\end{equation*}

Let
\[
   \B_{\operatorname{sp}}^+{}_a 
   \left(\text{resp.\ }\B_{\operatorname{sp}}^-{}_a\right)
   = 
   \left.\left\{
   T = 
\renewcommand{\arraycolsep}{0pt}
\begin{array}{|c|}
   \hline \hbox to \yh{\hfill$\scriptstyle i_1$\hfill} \\
   \hline \hbox to \yh{\hfill$\vdots$\hfill} \\
   \hline \hbox to \yh{\hfill$\scriptstyle{i_n}$\hfill} \\
   \hline
\end{array}_{a\ve^{1-n}}
   \right|
   \begin{gathered}
   i_p\in \mathbf B, \; i_1 \prec i_2 \prec \dots \prec i_n
   \\
   \text{$i$ and $\overline{i}$ do not appear simultaneously}
   \\
   \text{if $i_p = n$, $n-p$ is even (resp.\ odd)}
   \\
   \text{if $i_p = \overline{n}$, $n-p$ is odd (resp.\ even)}
   \end{gathered}
\right\}.
\]
We define $m_T$ as before. We have
\begin{equation}\label{eq:u2}
\begin{split}
   & u_{i,a\ve^{n-2p+i-1}}(m_T)
   = \left(i_p = i\right) - \left(i_{p+1} = i+1\right) \qquad
   \text{if $1\le i\le n-2$},
\\
   & u_{n-1,a\ve^{2n-2p-2}}(m_T) = \left(i_p = n-1\right)
     - \left(i_{p+1} = \overline{n}\right),
\\
   & u_{n,a\ve^{2n-2p-2}}(m_T) = - \left(i_p = \overline{n-1}\right)
     + \left(i_{p+1} = n \right).
\end{split}
\end{equation}

The proof of the following is left to the reader as an exercise.
\begin{Proposition}
\[
   \qch{q,t}(L(\Lambda_{n-1})_a)
    = \sum_{T\in \B_{\operatorname{sp}}^-{}_a} m_T,
    \qquad
   \qch{q,t}(L(\Lambda_{n})_a)
    = \sum_{T\in \B_{\operatorname{sp}}^+{}_a} m_T.
\]
\end{Proposition}

In particular, we find that $u_{i,a\ve^s}(m_T)$ is at most $1$ and if
$u_{i,a\ve^s}(m_T) = 1$, then $u_{i,b}(m_T) = 0$ for other $b$. This
implies that the graph $\Gamma_P$ is the same as the crystal graph.
More precisely, edges are given by
\begin{equation*}
\renewcommand{\arraycolsep}{0pt}
\begin{array}{|c|l}
    \hbox to \yh{\hfill\hfill} & \\
    \cline{1-1} \hbox to \yh{\hfill$\scriptstyle i$\hfill}
     & \;{}_{a\ve^{n+1-2p}}\\
    \cline{1-1} \hbox to \yh{\hfill$\vdots$\hfill} & \\
    \cline{1-1} \!\!\hbox to
    \yh{\hfill$\scriptstyle\overline{i\!+\!1}$\hfill}
        & \;{}_{a\ve^{n+1-2p'}} \\
    \cline{1-1} \hbox to \yh{\hfill\hfill} &
\end{array}
   \xrightarrow{i,a\ve^{n-2p+i}}
\renewcommand{\arraycolsep}{0pt}
\begin{array}{|c|}
   \hbox to \yh{\hfill\hfill} \\
   \hline \!\hbox to \yh{\hfill$\scriptstyle i\!+\!1$\hfill} \\
   \hline \hbox to \yh{\hfill$\vdots$\hfill} \\
   \hline \hbox to \yh{\hfill$\scriptstyle\overline{i}$\hfill} \\
   \hline \hbox to \yh{\hfill\hfill}
\end{array}
   \quad (1\le i\le n-1),
\qquad
\renewcommand{\arraycolsep}{0pt}
\begin{array}{|c|l}
   \hbox to \yh{\hfill\hfill} \\
   \cline{1-1} \!\!\hbox to \yh{\hfill$\scriptstyle n\!-\!1$\hfill}
   & \;{}_{a\ve^{n+1-2p}}\\
   \cline{1-1} \hbox to \yh{\hfill$\scriptstyle n$\hfill} & \\
   \cline{1-1} \hbox to \yh{\hfill\hfill} &
\end{array}
   \xrightarrow{n,a\ve^{n-1-2p}}
\renewcommand{\arraycolsep}{0pt}
\begin{array}{|c|}
   \hbox to \yh{\hfill\hfill} \\
   \hline \hbox to \yh{\hfill$\scriptstyle \overline{n}$\hfill} \\
   \hline \!\!\hbox to \yh{\hfill$\scriptstyle\overline{n\!-\!1}$\hfill} \\
   \hline \hbox to \yh{\hfill\hfill}
\end{array}
\;\;.
\end{equation*}

We have
\begin{equation}\label{eq:v2}
   v_{i,\ve^{n-2p+i}}(m_T, m_P) =
   \begin{cases}
     \left(p \le i \prec i_p\right) & \text{if $1\le i\le n-2$},
   \\
     \left(p \neq n \preceq i_p\right) & \text{if $i = n-1$},
   \\
     \left(\overline{p} \neq \overline{n} \preceq i_p\right)
     & \text{if $i = n$}.
   \end{cases}
\end{equation}

Suppose that a Drinfeld polynomial $P$ is given. We define the set of
column increasing tableaux $\B(P)$ as in the type $A_n$ case.
For a tableau $T = (T_1, T_2, \dots, T_L)\in\B(P)$, we define 
\(
  d(T_\alpha, T_\beta) \defeq
  d(m_{T_\alpha}, m_{P^\alpha}; m_{T_\beta}, m_{P^\beta})
\),
substituting (\ref{eq:ii}, \ref{eq:v1}, \ref{eq:u2},
\ref{eq:v2}) into \eqref{eq:d}. (We do not try to simplify the
expression as in the type $A_n$ case.)
Then we set $d(T) = \sum_{\alpha<\beta} d(T_\alpha, T_\beta)$ as
before. We also set $l(T) = \sum_\alpha l(T_\alpha)$, where
$l(T_\alpha)$ was defined in \eqref{eq:expo}. We get
\begin{Theorem}
\begin{equation*}
    \qch{q,t}(M(P))
    = \sum_{T\in\B(P)} t^{2d(T)+2l(T)} m_T.
\end{equation*}
\end{Theorem}

\end{document}